\documentclass[11pt]{article}

\usepackage{amsfonts}
\usepackage{amsfonts,amssymb,amsmath,amsbsy,amsthm,amscd}
\usepackage[utf8]{inputenc} 
\usepackage[english]{babel}   
\usepackage{indentfirst} 
\begin{document}

\begin{center}
RELIABILITY ANALYSIS OF THE SYSTEM LIFETIME
WITH $N$ WORKING ELEMENTS AND REPAIRING DEVICE UNDER CONDITION OF THE ELEMENT'S FAST REPAIR
\end{center}

\begin{center}
Golovastova E. A.
\end{center}

\textbf{Annotation.} In this paper we consider the system with $n$ identical elements and one repairing device. In each time moment only one element works while the rest elements stay in reserve. The distribution of element repairing period is exponential, the distribution of element working period is general. We obtain the asymptotic distribution of the system lifetime under condition of the element's fast repair.



\section{Introduction}
In this paper we consider the one of the problems from the mathematical reliability theory [1]. We suppose, that the system consists of $ n $ identical elements, one of which works, while the others stay in reserve. Also there is one secure repairing device. The device can restore only one element at an each moment of time, and if there are another broken elements, they stay in a queue. The element working time $\eta$ has general distribution, and its recovery time $\xi$ has an exponential distribution.  The system crashes at the time $\tau$, when all $n$ elements brakes down. All described random variables are supposed to be mutually independent.

The most important characteristics of this  model can be obtained by relying on the results from [1], [2]; the results for general similar models presented, for example, in the following papers [3], [4]. 

Since real technical systems nowadays are sufficiently reliable, there appears the problem  of an asymptotic analysis of systems with high reliability.

This paper established that in described model under the condition of element fast repairing: $ E  \xi  \to  \infty$ (so, this 
involves, that $P(\xi \leq t) \to 1$), and some $\varepsilon (E \xi) \to 0$, then the distribution of the random variable $ \varepsilon (E \xi)^{n-1} \cdot \tau $,  converges to exponential with parameter $1 \setminus E \eta$.

\section{Model description}
We consider a system consisting of $ n $ identical elements, where one element works, and the rest ones are in reserve. Working element can brake down and then it is immediately replaced by operable one. Repairing device is reliable. The system crashes at the moment when all $ n $ elements become inoperative. We will find equations for the distribution of system operating time and its asymptotic behaviour when the elements are quickly restored.

We assume, that all the random variables that determine system are mutually independent. Also, we suppose, that recovery times of the broken elements are exponentially distributed with the parameter $ \mu $, and elements working times have general distribution ~--~ $G(t)$. Let $G(t)$ be a continuous function for $t\geq 0$, and have finite moments. Particularly, let $\int_0^\infty t \ dG(t) = b < \infty$.

Denote by $\nu (\eta)$ ~--~ the number of broken elements, which have been repared during the element working time ~--~ $\eta$;
 $\tau_j$ ~--~ the system lifetime, when it starts working with $j$ broken elements, $j = 0,1,2,\ldots$.

Lets notice, that:
$$
P(\nu (\eta)=j)=e^{-\mu\eta}\frac{(\mu\eta)^j}{j!}, \quad j=0,1,2,\ldots
\eqno(1)
$$

We consider the system under the condition, that the elements are quickly restored. It means, that $ E  \xi = \mu  \to \infty$, or, in particular, it entails, that the following probability:
$$
P(\xi > \eta )=\int_0^\infty e^{-\mu t}\ dG(t) = \varepsilon(\mu) \to 0
\eqno(2)
$$
as $\mu  \to \infty$. Due to brevity, sometimes we will omit in $\varepsilon(\mu)$ the dependence on the argument below and write just $\varepsilon$.

We will show, that:
$$
P(\varepsilon(\mu)^{n-1} \tau_j > t) \to e^{-\frac{t}{b}} \text{ when } \mu \to \infty , \quad j=0,1,2,\ldots.
\eqno(3)
$$

By definition, $f(x)\sim g(x)$ when $x \to a$, if $\lim_{x \to a} \frac{f(x)}{g(x)} = 1$.

\section{System operating time distribution}
The system lifetime is determined by following stochastic equations (here $n >2$):
$$
\tau_0 = \eta + \tau_1
$$
$$
\tau_1 = (\eta + \tau_1)\ I(\nu(\eta) \geq 1) + (\eta + \tau_2)\ I(\nu(\eta) = 0)
$$
$$
\tau_j = (\eta + \tau_1)\ I(\nu(\eta) \geq j) + \sum_{k=0}^{j-1}(\eta + \tau_{j+1-k})I(\nu(\eta) = k), \quad 1< j < n-1
$$
$$
\tau_{n-1} =  (\eta + \tau_1)\ I(\nu(\eta) \geq n-1) + \sum_{k=1}^{n-2} (\eta + \tau_{n-k})I(\nu(\eta) = k) + \eta I(\nu(\eta) = 0)
$$
Further we take the  Laplace–Stieltjes transforms from previous equations and get the following linear algebraic system of equations:

$$
\varphi_0(s) = g(s)\varphi_1(s)
$$
$$
\varphi_1(s) = (g(s)-g_0(s))\varphi_1(s) + \varphi_2(s)g_0(s)
$$
$$
\varphi_j(s) = \left( g(s)-\sum_{k=0}^{j-1}g_k(s)\right) \varphi_1(s) + \sum_{k=0}^{j-1}g_k(s)\varphi_{j+1-k}(s), \quad 1< j < n-1
\eqno(4)
$$
$$
\varphi_{n-1}(s) = \left( g(s)-\sum_{k=0}^{n-2}g_k(s)\right) \varphi_1(s) + \sum_{k=1}^{n-2}g_k(s)\varphi_{n-k}(s) + g_0(s)
$$
where:
$$
\varphi_j(s) = E \ e^{-s\tau_j}, \quad g(s) = \int_0^\infty e^{-sx}\ dG(x), \quad g_0(s) = \int_0^\infty e^{-(s+\mu)x}\ dG(x)
$$
$$
g_j(s) = \int_0^\infty e^{-(s+\mu)x}\frac{(\mu x)^j}{j!}\ dG(x),\quad j = 1,2,\ldots .
$$
So, the convergence (3) is equivalent to convergence:
$$
\varphi_j(\varepsilon(\mu)^{n-1}s) \to \frac{1}{1+bs} \text{ when } \ \mu \to \infty, \quad j=0,1,2,\ldots.
\eqno(5)
$$
Taking into account (2), we notice, that when $\mu \to \infty$:
$$
g(\varepsilon^{n-1} s) \sim 1 - \varepsilon ^{n-1}bs
$$
$$
g_0(\varepsilon^{n-1} s) \sim \varepsilon - \varepsilon ^{n-1}s \int_0^\infty x e^{-\mu x}\ dG(x) = \varepsilon - \varepsilon ^{n-1}s \ \gamma_0(\mu)
\eqno(6) 
$$
$$
g_j(\varepsilon^{n-1} s) \sim  g_j(0) - \varepsilon ^{n-1}s \int_0^\infty x\frac{(\mu x)^j}{j!} e^{-\mu x}\ dG(x) = g_j(0) - \varepsilon ^{n-1}s \ \gamma_j(\mu)
\quad j=1,2,\ldots
$$
we also mention, that:
$$
g_j(0)=g_j(\mu , 0) \to 0, \ j=1,2,\ldots 
$$
$$
\gamma_j=\gamma_j(\mu) \to 0, \ j=0,1,2,\ldots
$$
when $\mu \to \infty$.

\section{Asymptotic analysis of the operating time assuming its high reliability}
For $2$ elements in the system, we have the following equations for  Laplace–Stieltjes transform of system (4):
$$
\varphi_0(s) =  g(s) \varphi_1(s)
$$
$$
\varphi_1(s) = (g(s)-g_0(s))\varphi_1(s) + g_0(s).
$$
$$
\varphi_1(s) = \frac{g_0(s)}{1-g(s)+g_0(s)}, \quad \varphi_0(s) =  g(s) \varphi_1(s).
$$
When $\mu \to \infty$ from (6):
$$
g(\varepsilon s) \sim 1 - \varepsilon bs
$$
$$
g_0(\varepsilon s) \sim \varepsilon - \varepsilon s \ \gamma_0 \sim \varepsilon.
$$
So, when $\mu \to \infty$ functions $\varphi_0(\varepsilon s), \varphi_1(\varepsilon s) $ tend to $1\setminus(1+bs)$, and (5) is true for $2$ elements in the system. Further we assume, that $n>2$.
When $\mu \to \infty$:
$$
\frac{\varphi_0(\varepsilon^{n-1}s)}{\varphi_1(\varepsilon^{n-1}s)}=g(\varepsilon^{n-1}s)
$$
$$
\frac{\varphi_0(\varepsilon^{n-1}s)}{\varphi_1(\varepsilon^{n-1}s)} \sim 1-\varepsilon^{n-1}sb.
$$
So, $\varphi_0(\varepsilon^{n-1}s) \sim \varphi_1(\varepsilon^{n-1}s)$. Due to brevity, we will sometimes omit the dependence of the argument $\varepsilon^{n-1}s$ below.
$$
1=(g-g_0) + \frac{\varphi_2}{\varphi_1}g_0
$$
$$
1 \sim 1 - \varepsilon^{n-1}sb - \varepsilon + \varepsilon^{n-1}s \gamma_0 + \frac{\varphi_2}{\varphi_1}\cdot \left(\varepsilon - \varepsilon^{n-1}s \gamma_0 \right)
$$
Leaving the slowest decreasing terms with $\varepsilon$ in the last expression, we get:
$$
\frac{\varphi_2}{\varphi_1} \sim
1+\varepsilon^{n-2}sb
$$
So, $\varphi_2 \sim \varphi_1$ and $\varphi_2 \sim \varphi_1 \cdot (1+\varepsilon^{n-2}sb)$. For $2<j<n-1$ we get:
$$
\frac{\varphi_j}{\varphi_1}=\left( g - g_0 - g_1 - \ldots - g_{j-1} \right) +
$$
$$
+\frac{\varphi_2}{\varphi_1}g_{j-1}+\frac{\varphi_2}{\varphi_1}g_{j-2}+\ldots
+\frac{\varphi_j}{\varphi_1}g_{1}+\frac{\varphi_{j+1}}{\varphi_1}g_{0}
$$
Again, using (6), we leave the slowest decreasing terms with $\varepsilon$.
$$
1+\varepsilon^{n-j}sb \sim 1 - \varepsilon + g_1(0)\cdot \varepsilon^{n-j}sb + \frac{\varphi_{j+1}}{\varphi_1} \cdot (\varepsilon - \varepsilon^{n-1}s \gamma_0)
$$
$$
\frac{\varphi_{j+1}}{\varphi_1} \sim
1+\varepsilon^{n-j-1}sb
$$
So, $\varphi_{j+1} \sim \varphi_1$ and $\varphi_{j+1} \sim \varphi_1 \cdot (1+\varepsilon^{n-j-1}sb)$. Particularly, $\varphi_{n-1} \sim \varphi_1$ and $\varphi_{n-1} \sim \varphi_1 \cdot (1+\varepsilon \cdot sb)$. So, we got, that all $\varphi_{k}(\varepsilon^{n-1}s) \sim \varphi_1 (\varepsilon^{n-1}s) $, $k=2,3,\cdots , n-1$, when $\mu \to \infty$. Using this result, we divide the last equation in system (4) by $\varphi_1$, and, leaving the slowest decreasing terms with $\varepsilon$, we get:
$$
\varphi_1 \cdot(1 + \varepsilon \cdot sb) \sim \varphi_1 \cdot (1-\varepsilon + \varepsilon  sb \cdot g_1(0)) + (\varepsilon - \varepsilon^{n-1}s \gamma_0)
$$
$$
\varphi_1 \cdot(\varepsilon + \varepsilon \cdot sb - \varepsilon \cdot sb \cdot g_1(0) ) \sim \varepsilon \cdot(1-\varepsilon^{n-2}s \gamma_0)
$$
So, $\varphi_1(\varepsilon^{n-1}s) \sim 1\setminus(1+bs)$, when $\mu \to \infty$, and, considering all above, (5) again is true.
\\
\\
\begin{center}
\textbf{Literature}
\end{center}

1. Gnedenko B.V., Belyaev Yu.K., Soloviev A.D.
Mathematical methods in the theory of reliability. M.: Nauka, 1965.

2. Hutson V.C.L, Pym J.S.
Applications of Functional Analysis and Operator Theory. Academic Press, London, New York, Toronto, 1980, 432p. 

3. Golovastova Eleonora A. System operating time with unreliable interchangeable elements and precarious recovery device //
Tomsk State University Journal of Control and Computer Science. 2020. N 52. pp. 59-65.

4. Golovastova E.A. The system operating time with two different unreliable servicing devices // 2018. URL: https://arxiv.org/abs/1811.12193

\end{document}